\documentclass{amsart}
\usepackage{amssymb}
\newtheorem{Proposition}{Proposition}
\newtheorem{Corollary}{Corollary}
\newtheorem{Theorem}{Theorem}
\theoremstyle{definition}
\newtheorem*{ack}{Acknowledgment}
\theoremstyle{remark}
\newtheorem{Conjecture}{Conjecture}
\numberwithin{equation}{section}

\DeclareMathOperator{\End}{End}
\DeclareMathOperator{\Tr}{Tr}
\DeclareMathOperator{\id}{id}
\def\A{\mathsf{A}}
\def\B{\mathsf{B}}
\def\Bb{\tilde{\mathsf{B}}}

\def\E{\mathsf{E}}
\def\F{\mathsf{F}}
\def\G{\mathsf{G}}
\def\H{\mathsf{H}}
\def\GG{\mathbf{G}}
\def\K{\mathsf{K}}
\def\L{\mathsf{L}}
\def\M{\mathsf{M}}
\def\X{\mathsf{X}}

\def\Y{\mathsf{Y}}
\def\C{\mathbb{C}}
\def\S{\mathsf{S}}
\def\T{\mathsf{T}}
\def\P{\mathsf{P}}
\def\u{\mathsf{u}}
\def\v{\mathsf{v}}
\def\w{\mathsf{w}}
\def\pw{\tilde{\mathsf{w}}}
\def\q{\mathsf{q}}
\def\Z{\mathsf{Z}}
\def\x{\mathsf{x}}
\def\y{\mathsf{y}}
\def\z{\mathsf{z}}
\def\xx{\mathbf{x}}
\def\yy{\mathbf{y}}
\def\e{\mathsf{e}}
\def\a{\mathsf{a}}
\def\f{\mathsf{f}}
\def\RR{\mathbf{R}}
\def\Rr{\mathcal{R}}
\def\y{\mathsf{y}}
\def\n{\mathsf{N}}
\def\b{\mathsf{b}}
\def\R{\mathsf{R}}
\def\U{\mathsf{U}}
\def\V{\mathsf{V}}
\def\W{\mathsf{W}}
\def\Uu{{\mathcal U}}

\title{From the tetrahedron equation to universal $R$-matrices}
\author{R.~M.~Kashaev}
\address{Helsinki Institute of Physics}
\curraddr{St.~Petersburg Steklov Mathematical Institute}
\email{kashaev@pdmi.ras.ru}
\author{A.~Yu.~Volkov}
\address{Helsinki Institute of Physics}
\email{volkov@rock.helsinki.fi}
\date{December 1998}
\keywords{Tetrahedron equation, quantum group, universal $R$-matrix}

\begin{document}
\begin{abstract}
Modified universal $R$-matrices,
associated with the central extension
(through the Drinfeld's double construction) of the quantum groups
$\U_\q(\mathfrak{sl}_\n)$,
are realized through an infinite dimensional spectral parameter
dependent solution for the tetrahedron equation, provided a certain
identity on $q$-exponentials holds true.
\end{abstract}

\maketitle

\section{Introduction and notation}

The (quantum) Yang--Baxter equation
(YBE) \cite{Yang,Baxter} is well known to play the fundamental role
in constructing and solving integrable models of two-dimensional
statistical physics and quantum field theory \cite{Baxter1,Faddeev}.
The tetrahedron equation (TE)\cite{Zamolodchikov} has been introduced
as a three-dimensional analogue of the YBE. One of the features of the TE
is the possibility to construct from one solution of it an
infinite sequence of solutions of the YBE. This means that
one integrable three-dimensional lattice model combines an
infinite family of integrable two-dimensional models. A notable example
is the Zamolodchikov--Baxter--Bazhanov (ZBB) three-dimensional
lattice model \cite{Zamolodchikov,BB1,BB2}:
on a cubic lattice, with one direction being periodic of finite length $\n$,
this model is equivalent to the two-dimensional
$\mathfrak{sl}_\n$
 chiral Potts model \cite{BKMS}. The field theoretical
counterpart of the ZBB model is the
three-dimensional bilinear Hirota equation
\cite{Hirota} which, being considered on a lattice
periodic in one direction with period $\n$, can be interpreted
as a discrete  version of the affine
$\mathfrak{sl}_\n$ Toda field theory. The quantum theory of
this system has been developed in \cite{FV1,FV2}
and in \cite{Bobenko}
for the case $\n=2$, and subsequently in \cite{KR} for any $\n$.
Thus, the TE can be a powerful unifying tool for different
integrable models in two dimensions.

In \cite{Sergeevetal} few nontrivial examples of formal infinite
dimensional solutions for the TE have been found.
In the classical limit, when the deformation parameter $\q$ tends to
unity, these solutions are related to
(infinite dimensional) functional solutions for the TE
\cite{Kashaev,Sergeev1,KKS}, which
in turn are related to solutions for the local YBE
\cite{MN,Maillet} --- the three-dimensional counterpart
of the zero-curvature condition.
Using one of the solutions from \cite{Sergeevetal},
Sergeev in \cite{Sergeev} obtained a family of
(spectral parameter dependent) formal
infinite dimensional solutions for the YBE and
interpreted them as specializations of the affine universal $R$-matrices.

The purpose of this paper is to uncover the three-dimensional nature
of the quantum groups $\U_\q(\mathfrak{sl}_\n)$
as well as the corresponding
universal $R$-matrices.
First, we rigorously define a close (spectral parameter dependent)
cousin of one of the solutions of \cite{Sergeevetal} as a formal
power series in the spectral parameter with operator coefficients
acting in the space of Laurent
polynomials of three indeterminates.
We use this solution to compose a solution of the YBE
as a three dimensional analogue of the transfer matrix
of the size $\n\times\n\times\n$. Then, assuming validity of a certain
identity on $q$-exponentials, we identify this with the image
of the centrally extended $\U_\q(\mathfrak{sl}_\n)$
(the Drinfeld double of the Borel subalgebra)
modified universal $R$-matrix under an infinite dimensional representation
in the space of Laurent polynomials of $\n^2$ indeterminates. This
representation is conjectured to be faithful.
A finite dimensional variant of such representation in the
case of roots of unity first has been constructed by Tarasov
in \cite{Tarasov}. Restriction of this representation
to the Borel subalgebra is known as Feigin's homomorphism \cite{Ber} or
free field representation \cite{Mor}.
The universal $R$-matrix itself is not well defined
in this representation, since
the power series expansion of the $q$-exponentials, which enter
the structure of the universal $R$-matrix, do not truncate. We modify
the definition of the universal $R$-matrix by introducing $\n-1$
fictitious spectral parameters and consider it
as a generating function for operator
coefficients of the power series expansion in these spectral parameters.
It is this modified universal $R$-matrix the image of which we identify
with our three dimensional construction.

 We fix notation.
 Let $\q$ be an indeterminate, and
let $\A\equiv\End(\C(\q)[t,t^{-1}])$
be the algebra of linear operators in
the space of Laurent polynomials in one
indeterminate with coefficients from the field $\C(\q)$.
 By abuse of notation, we denote
        \[
\A^{\otimes n}\equiv \End(\C(\q)[t,t^{-1}]^{\otimes n})\equiv
\End(\C(\q)[t_1,t_1^{-1},\ldots,t_n,t_n^{-1}]),
        \]
and define similarly algebras $\A^{\otimes m}\otimes\A^{\otimes n}$ and
$(\A^{\otimes m})^{\otimes n}$.

Let $\B$ be any algebra. For any positive $m$ define
natural embeddings
        \[
\iota_i\colon\B\ni b\mapsto 1\otimes\cdots\otimes
b\otimes\cdots\otimes1\in\B^{\otimes m},
        \]
where element $b$ in the
r.h.s. stands on the $i$-th position.  For a finite sequence of such
embeddings
        \[
\B\stackrel{\iota_{i_1}}{\longrightarrow}\B^{\otimes
n_1} \stackrel{\iota_{i_2}}{\longrightarrow} (\B^{\otimes
n_1})^{\otimes n_2} \stackrel{\iota_{i_3}}{\longrightarrow}\cdots
\stackrel{\iota_{i_s}}{\longrightarrow}(\ldots(\B^{\otimes
n_1})\ldots)^{\otimes n_s}
        \]
 we denote by
        \[
b_{i_s:\ldots:
i_2:i_1}\equiv\iota_{i_s:\ldots:i_2:i_1}(b)\equiv
\iota_{i_s}\circ\cdots\circ \iota_{i_2}
\circ\iota_{i_1}(b)
        \]
the image of element $b\in\B$.
  In the
notation for these embeddings we intentionally suppress any indication
of the algebras involved since these will be clear in each concrete
case from the context.

More generally, if $\B_1,\ldots,\B_l$ are algebras, we shall write
        \[
b_{i_{11}:\ldots:i_{1s_1},\ldots,i_{l1}:\ldots:i_{ls_l}}\equiv
\iota_{i_{11}:\ldots:i_{1s_1}}\otimes\cdots\otimes
\iota_{i_{l1}:\ldots:i_{ls_l}}(b)
        \]
for any $b\in\B_1\otimes\cdots\otimes\B_l$.

We shall encounter products of noncommuting elements, so to fix
the order we shall write
        \[
  \prod_{i\uparrow_1^n}\alpha_i\equiv\alpha_1\alpha_2\cdots\alpha_n,
\quad
\prod_{i\downarrow_1^n}\alpha_i\equiv\alpha_n\alpha_{n-1}\cdots\alpha_1.
        \]

\section{A solution for the (spectral parameter dependent)
tetrahedron equation}

Define elements $\u$, $\v\in \A$:
        \[
\u(f(t))\equiv f(\q t),\quad
 \v(f(t))\equiv tf(t),\quad \forall f(t)\in \C(\q)[t,t^{-1}],
        \]
 and
element $\F\in\A^{\otimes 3}$:
        \begin{equation}\label{Fdefinition}
 \F(f(t_1,t_2,t_3))\equiv
f(t_1t_2/t_3,t_3,t_2),\quad \forall f(t_1,t_2,t_3)\in
(\C(\q)[t,t^{-1}])^{\otimes 3}
        \end{equation}
        \begin{Proposition}
Operators $\u$,
$\v$, and $\F$ satisfy the following relations:
        \begin{gather}
\u\v=\q\u\v,\quad\F^2=1,\notag\\
       \label{uvFcommutation}
\u_1\F=\F u_1,\quad\u_2\F=\F\u_1\u_3,\quad
\v_1\v_2\F=\F\v_1\v_2,\quad\v_2\F=\F\v_3,\\
        \label{Ftetrahedron}
\F_{1,2,4}\F_{1,3,5}\F_{2,3,6}\F_{4,5,6}=\F_{4,5,6}\F_{2,3,6}\F_{1,3,5}
\F_{1,2,4}.
        \end{gather}
        \end{Proposition}
\begin{proof}
It is a straightforward verification.
\end{proof}

Eqn~(\ref{Ftetrahedron}) is called the (constant)
TE and the
operator $\F$ defined by eqn~(\ref{Fdefinition}) is one of the
simplest examples of functional solutions to the TE
see \cite{Kashaev,Sergeev1,KKS}.
Define now an invertible element
$\B^\x\in\A^{\otimes 3}[[\x]]$ in the algebra of formal power series
 with coefficients from $\A^{\otimes 3}$:
        \begin{equation}\label{Bdefinition}
\B^\x=\F\psi(\x\v_1\u^{-1}_1\u_2
 \u^{-1}_3\v^{-1}_3)=\sum_{k=0}^\infty \x^k \b_k,
        \end{equation}
 where
        \[
\psi(\x)\equiv (-\x;\q^2)_\infty^{-1}=\sum_{k=0}^\infty
\frac{(-\x)^k}{(\q^2;\q^2)_k}, \quad
 \b_k=\F\frac{(-\v_1\u^{-1}_1
\u_2\u^{-1}_3\v^{-1}_3)^k}{(\q^2;\q^2)_k}\in \A^{\otimes 3},
        \]
with the standard notation in $q$-mathematics
        \[
(\x;\y)_k\equiv\prod_{i=0}^{k-1}(1-\x\y^i).
        \]
One can look at $\B^\x$ as the generating function for the infinite
sequence of operators $\b_k$.

We shall use the following five-term identity satisfied by the
$\psi$-function:
        \begin{equation}\label{pentagon}
 \psi(\X)\psi(\Y)=\psi(\Y)\psi(\Y\X)\psi(\X)
        \end{equation}
where operators $\X$ and $\Y$ are such that
$\X\Y=\q^2\Y\X$. Note that this identity is understood as an identity
for formal power series in noncommuting quantities $\X$ and $\Y$.
        \begin{Proposition}
The following spectral parameter dependent
TE in $\A^{\otimes 6}[[x,y,z]]$ holds
        \begin{equation}\label{Btetrahedron}
\B_{1,2,4}^\x\B_{1,3,5}^{\x\z}\B_{2,3,6}^\y \B_{4,5,6}^\z=
\B_{4,5,6}^\y\B_{2,3,6}^\z\B_{1,3,5}^{\x\y}\B_{1,2,4}^\x.
        \end{equation}
        \end{Proposition}
\begin{proof}
Substituting definition~(\ref{Bdefinition}) into
eqn~(\ref{Btetrahedron}), one can cancel all the
$\F$-operators by first moving them to the left, using
eqns~(\ref{uvFcommutation}), and then applying the
TE~(\ref{Ftetrahedron}). What remains is the
following operator identity:
        \begin{equation}\label{fourterm}
\psi(\U)\psi(\U\W)\psi(\V)\psi(\W)=\psi(\V)\psi(\W)\psi(\V\U)\psi(\U),
        \end{equation}
where operator combinations
        \[
\U=x\v_1\u_1^{-1}\u_2\u_3^{-1}\v_3^{-1},\quad
\V=y\v_2\u_2^{-1}\u_3\u_4\u_5^{-1}\v_5^{-1},\quad
\W=z \v_3\u_3^{-1}\u_5\u_6^{-1}\v_6^{-1}
        \]
satisfy relations
        \[
\U\V=\q^2\V\U,\quad  \V\W=\q^2\W\V,\quad \W\U=\q^2\U\W.
        \]
The latter relations together with the five-term
identity (\ref{pentagon}) imply eqn~(\ref{fourterm})
through the following sequence of equalities (in each step
the fragment to be transformed is underlined):
        \begin{multline*}
\psi(\U)\psi(\U\W)\underline{\psi(\V)\psi(\W)}=
\underline{\psi(\U)\psi(\U\W)\psi(\W)}\psi(\W\V)\psi(\V)\\
=\psi(\W)\underline{\psi(\U)\psi(\W\V)}\psi(\V)=
\psi(\W)\psi(\W\V)\underline{\psi(\U)\psi(\V)}\\
=\underline{\psi(\W)\psi(\W\V)\psi(\V)}\psi(\V\U)\psi(\U)=
\psi(\V)\psi(\W)\psi(\V\U)\psi(\U).
        \end{multline*}
\end{proof}

Formally, at the special value of the spectral parameter $\x=1$
we obtain solution for the constant TE.
Slightly different form of this solution first has been found
in \cite{Sergeevetal}.

Note another special value of the spectral parameter
when we get (trivial) constant solution
$\B^{\x=0}=\F.$

\section{Solution for the
Yang--Baxter equation}

Fix a positive integer $\n\ge2$
and define
$\Bb^\xx\in((\A^{\otimes\n})^{\otimes2}\otimes\A)[[\xx]]$,
where $\xx\equiv(\x_1,\ldots,\x_\n)$,
        \begin{equation}\label{kulbasa}
\Bb^\xx=\prod_{i\downarrow_1^\n}\B_{1:i,2:i,3}^{\x_i}.
        \end{equation}
In the rest of the paper all indices, taking values $1,\ldots,\n$,
will be considered $\pmod\n$.
        \begin{Proposition}
The following equation is satisfied  in $((\A^{\otimes\n})^{\otimes3}
\otimes\A^{\otimes3})[[\xx,\yy]]$
        \begin{equation}\label{kulbasaYBE}
\Bb_{1,2,4}^\xx\Bb_{1,3,5}^{\xx\tau(\yy)}
\Bb_{2,3,6}^\yy
\B_{4,5,6}^{\y_\n}=\B_{4,5,6}^{\y_\n}
\Bb_{2,3,6}^{\tau(\yy)}\Bb_{1,3,5}^{\xx\yy}\Bb_{1,2,4}^\xx,
        \end{equation}
where
$\xx\yy\equiv(\x_1\y_1,\ldots,\x_\n\y_\n)$
and $\tau(\xx)\equiv(\x_\n,\x_1,\ldots,\x_{\n-2},\x_{\n-1})$.
        \end{Proposition}
\begin{proof} Substituting definition~(\ref{kulbasa}) into the l.h.s.
 of eqn~(\ref{kulbasaYBE}), we have the following sequence
of equalities (the fragments eqn~(\ref{Btetrahedron}) to be applied to
are underlined):
\begin{multline*}
\Bb_{1,2,4}^\xx\Bb_{1,3,5}^{\xx\tau(\yy)}\Bb_{2,3,6}^\yy\B_{4,5,6}^{y_\n}=
\prod_{i\downarrow_1^\n}\left(\B_{1:i,2:i,4}^{\x_i}
\B_{1:i,3:i,5}^{\x_i\y_{i-1}}
\B_{2:i,3:i,6}^{\y_i}\right) \B_{4,5,6}^{y_\n}\\
=\prod_{i\downarrow_2^{\n}}
\left(\B_{1:i,2:i,4}^{\x_i}
\B_{1:i,3:i,5}^{\x_i\y_{i-1}}
\B_{2:i,3:i,6}^{\y_i}\right)
\underline{\B_{1:1,2:1,4}^{\x_1}\B_{1:1,3:1,5}^{\x_1\y_\n}
 \B_{2:1,3:1,6}^{\y_1} \B_{4,5,6}^{\y_\n}}\\
=\prod_{i\downarrow_3^\n}\left(
\B_{1:i,2:i,4}^{\x_i}
\B_{1:i,3:i,5}^{\x_i\y_{i-1}}
\B_{2:i,3:i,6}^{\y_i}\right)\\
\times
\underline{\B_{1:2,2:2,4}^{\x_2}
\B_{1:2,3:2,5}^{\x_2\y_1}
\B_{2:2,3:2,6}^{\y_2} \B_{4,5,6}^{\y_1}}
\B_{2:1,3:1,6}^{\y_\n}\B_{1:1,3:1,5}^{\x_1\y_1}
\B_{1:1,2:1,4}^{\x_1}
=\ldots\\=\B_{4,5,6}^{y_\n}\prod_{i\downarrow_1^\n}
\left(\B_{2:i,3:i,6}^{\y_{i-1}}
\B_{1:i,3:i,5}^{\x_i\y_i}
\B_{1:i,2:i,4}^{\x_i}\right)
=\B_{4,5,6}^{y_\n}\Bb_{2,3,6}^{\tau(\yy)}\Bb_{1,3,5}^{\xx\yy}\Bb_{1,2,4}^\xx
\end{multline*}
thus obtaining the r.h.s. of eqn~(\ref{kulbasaYBE}).
\end{proof}

Operator $\Bb^\xx$ can be considered as a three-dimensional analogue
of the monodromy operator. The operator
\begin{equation}\label{Rdefinition}
\R^\xx\equiv\Tr_3\Bb^\xx\in (\A^{\otimes \n})^{\otimes2}[[\xx]],
\end{equation}
provided the trace does exist (see below), is then analogous
to the transfer-matrix.
\begin{Corollary}
Operator $\R^\xx$ satisfies the following YBE
\[
 \R_{1,2}^\xx\R_{1,3}^{\xx\tau(\yy)}\R_{2,3}^\yy
=\R_{2,3}^{\tau(\yy)}\R_{1,3}^{\xx\yy}\R_{1,2}^\xx.
\]
\end{Corollary}
This is multi spectral parameter YBE
with unusual difference property.
The usual difference property is achieved for the
combination
\[
\RR^\xx\equiv
\prod_{i\downarrow_1^\n,j\uparrow_1^\n}\R_{1:i,2:j}^{\tau^i(\xx)}
\in ((\A^{\otimes \n})^{\otimes \n})^{\otimes2}[[\xx]].
\]
 \begin{Proposition}
The following YBE is satisfied
\begin{equation}\label{YBE}
\RR_{1,2}^\xx\RR_{1,3}^{\xx\yy}\RR_{2,3}^\yy=
\RR_{2,3}^\yy\RR_{1,3}^{\xx\yy}\RR_{1,2}^\xx
\end{equation}
\end{Proposition}
\begin{proof}
Operator $\RR^\xx$ can be represented as
\begin{equation}\label{Lfactorization}
\RR^\xx=\prod_{j\uparrow_1^\n}\M_{1,2:j}^\xx,\quad
 \M^\xx\equiv \prod_{i\downarrow_1^\n}\R_{1:i,2}^{\tau^i(\xx)}\in
((\A^{\otimes\n})^{\otimes\n})\otimes(\A^{\otimes\n})[[\xx]].
\end{equation}
We have
\begin{multline*}
\M_{1,2}^\xx\M_{1,3}^{\xx\yy}\R_{2,3}^\yy=
\prod_{i\downarrow_1^\n}\left(\R_{1:i,2}^{\tau^i(\xx)}
\R_{1:i,3}^{\tau^i(\xx\yy)}\right)\R_{2,3}^\yy\\
=\prod_{i\downarrow_2^\n}\left(\R_{1:i,2}^{\tau^i(\xx)}
\R_{1:i,3}^{\tau^i(\xx\yy)}\right)
\underline{\R_{1:1,2}^{\tau(\xx)}
\R_{1:1,3}^{\tau(\xx\yy)}\R_{2,3}^\yy}\\
=
\prod_{i\downarrow_3^\n}\left(\R_{1:i,2}^{\tau^i(\xx)}
\R_{1:i,3}^{\tau^i(\xx\yy)}\right)
\underline{\R_{1:2,2}^{\tau^2(\xx)}
\R_{1:2,3}^{\tau^2(\xx\yy)}\R_{2,3}^{\tau(\yy)}}
\R_{1:1,3}^{\tau(\xx)\yy}\R_{1:1,2}^{\tau(\xx)}=\ldots\\
=\R_{2,3}^\yy\prod_{i\downarrow_1^\n}\left(
\R_{1:i,3}^{\tau^i(\xx\tau^{-1}(\yy))}
\R_{1:i,2}^{\tau^i(\xx)}\right)=
\R_{2,3}^\yy\M_{1,3}^{\xx\tau^{-1}(\yy)} \M_{1,2}^\xx.
 \end{multline*}
Using this identity together with
formula~(\ref{Lfactorization}), we obtain
\begin{multline*}
\RR_{1,2}^\xx\M_{1,3}^{\xx\yy}\M_{2,3}^\yy=
\prod_{i\uparrow_1^\n}\M_{1,2:i}^\xx\M_{1,3}^{\xx\yy}
\prod_{j\downarrow_1^\n}\R_{2:j,3}^{\tau^j(\yy)}\\
=\prod_{i\uparrow_1^{\n-1}}\M_{1,2:i}^\xx
\underline{\M_{1,2:\n}^\xx\M_{1,3}^{\xx\yy}
\R_{2:\n,3}^\yy}\prod_{j\downarrow_1^{\n-1}}\R_{2:j,3}^{\tau^j(\yy)}\\
=\R_{2:\n,3}^\yy\prod_{i\uparrow_1^{\n-2}}\M_{1,2:i}^\xx
\underline{\M_{1,2:\n-1}^\xx\M_{1,3}^{\xx\tau^{-1}(\yy)}
\R_{2:\n-1,3}^{\tau^{-1}(\yy)}}\prod_{j\downarrow_1^{\n-2}}
\R_{2:j,3}^{\tau^j(\yy)}\M_{1,2:\n}^\xx=\ldots\\
=\prod_{j\downarrow_1^\n}\R_{2:j,3}^{\tau^j(\yy)}\M_{1,3}^{\xx\yy}
\prod_{i\uparrow_1^\n}\M_{1,2:i}^\xx=
\M_{2,3}^\yy\M_{1,3}^{\xx\yy}\RR_{1,2}^\xx,
\end{multline*}
which implies eqn~(\ref{YBE}):
\begin{multline*}
\RR_{1,2}^\xx\RR_{1,3}^{\xx\yy}\RR_{2,3}^\yy=
\RR_{1,2}^\xx\prod_{i\uparrow_1^\n}\left(\M_{1,3:i}^{\xx\yy}
\M_{2,3:i}^\yy\right)
=
\underline{\RR_{1,2}^\xx\M_{1,3:1}^{\xx\yy}\M_{2,3:1}^\yy}
 \prod_{i\uparrow_2^\n}\left(\M_{1,3:i}^{\xx\yy}
\M_{2,3:i}^\yy\right)\\
=
\M_{2,3:1}^\yy\M_{1,3:1}^{\xx\yy}
\underline{\RR_{1,2}^\xx\M_{1,3:2}^{\xx\yy}\M_{2,3:2}^\yy}
 \prod_{i\uparrow_3^\n}\left(\M_{1,3:i}^{\xx\yy}
\M_{2,3:i}^\yy\right)=\ldots\\
=
\prod_{i\uparrow_1^\n}\left(\M_{2,3:i}^\yy
\M_{1,3:i}^{\xx\yy}\right)\RR_{1,2}^\xx=
\RR_{2,3}^\yy\RR_{1,3}^{\xx\yy}\RR_{1,2}^\xx.
\end{multline*}
\end{proof}

\section{Connection with a universal $R$-matrix}
We evaluate explicitly the
case given by the following choice of the multi spectral parameter
$\xx=(\x_1,\x_2,\ldots,\x_{\n-1},0)$.

\subsection{Calculation of the operator $\R^\xx$}
We shall find it convenient to denote
        \[
\w\equiv\v\u,\quad \pw\equiv\v\u^{-1}.
        \]
        \begin{Proposition}
If $\xx=(\x_1,\ldots,\x_{\n-1},0)$, the operator $\R^\xx$
defined by (\ref{Rdefinition}) has the following explicit form
        \begin{equation}\label{Rcalculated}
\R^{(\x_1,\ldots,\x_{\n-1},0)}=\G\prod_{i\downarrow_1^{\n-1}}
\psi(\x_i\pw_{1:i}\u_{1:i-1}\u_{2:i}\w_{2:i-1}^{-1}),
        \end{equation}
where operator $\G\in (\A^{\otimes\n})^{\otimes2}$ is defined by
        \begin{equation}\label{Hdefinition}
H(f(\ldots,t_{1:i},\ldots,t_{2:j},\ldots))=
f(\ldots,t_{1:i}t_{2:i}/t_{2:i+1},\ldots,t_{2:j+1},\ldots)
        \end{equation}
        \end{Proposition}
\begin{proof}
First, note that the operator $\F$ defined in eqn~(\ref{Fdefinition})
can be factorized in the form:
\[
\F=\S_{1,3}^{-1}\P_{2,3}\S_{1,3},
\]
where
$\S,\P\in\A^{\otimes2}$ are defined by
\[
\S(f(t_1,t_2))=f(t_1t_2,t_2),\quad \P(f(t_1,t_2))=f(t_2,t_1).
\]
Operator $\S$ satisfies relations:
        \begin{gather*}
\S_{1,2}\S_{1,3}=\S_{1,3}\S_{1,2} ,\quad
\S_{1,3}\S_{2,3}= \S_{2,3} \S_{1,3}, \\
\S \u_1=\u_1\S,\quad \S \v_1=\v_1\v_2\S,\quad \S \u_1\u_2=\u_2\S,\quad
\S \v_2=\v_2\S,
        \end{gather*}
and we have the following identity for the
partial trace of the operator $\P$:  \[  \Tr_2\P=1_\A\in\A.  \]
Now we calculate
 \begin{multline*} \R^{(\x_1,\ldots,\x_{\n-1},0)}=
\Tr_3\S_{1:\n,3}^{-1}\P_{2:\n,3}\S_{1:\n,3}\prod_{i\downarrow_1^{\n-1}}
\B_{1:i,2:i,3}^{\x_i}\\
=\Tr_3\S_{1:\n,3}^{-1}\S_{1:\n,2:\n}\prod_{i\downarrow_1^{\n-1}}
\B_{1:i,2:i,2:\n}^{\x_i}\P_{2:\n,3}
=\S_{1:\n,2:\n}\prod_{i\downarrow_1^{\n-1}}
\B_{1:i,2:i,2:\n}^{\x_i}\S_{1:\n,2:\n}^{-1}\\
=\S_{1:\n,2:\n}\prod_{i\downarrow_1^{\n-1}}\F_{1:i,2:i,2:\n}\\
 \times
\prod_{j\downarrow_2^{\n-1}}\left(\prod_{k\uparrow_1^{j-1}}
\F_{1:k,2:k,2:\n}^{-1}
\psi(\x_j\pw_{1:j}\u_{2:j}\w_{2:\n}^{-1})
\prod_{l\downarrow_1^{j-1}}
\F_{1:l,2:l,2:\n}\right)\\
\times\psi(\x_1\pw_{1:1}\u_{2:1}\w_{2:\n}^{-1})
\S_{1:\n,2:\n}^{-1}\\
=\S_{1:\n,2:\n}\prod_{i\downarrow_1^{\n-1}}\F_{1:i,2:i,2:\n}
\S_{1:\n,2:\n}^{-1}
\prod_{j\downarrow_1^{\n-1}}
\psi(\x_j\pw_{1:j}\u_{1:j-1}\u_{2:j}\w_{2:j-1}^{-1}).
\end{multline*}
Thus, we come to formula~(\ref{Rcalculated}), where
\[
\G=\S_{1:\n,2:\n}\prod_{i\downarrow_1^{\n-1}}\F_{1:i,2:i,2:\n}
\S_{1:\n,2:\n}^{-1} =\prod_{i\downarrow_1^{\n-1}}\P_{2:i,2:\n}
\prod_{j=1}^\n\left( \S_{1:j,2:j}^{-1}\S_{1:j,2:j-1}\right),
\]
and it is easy to check up eqn~(\ref{Hdefinition}).
\end{proof}

Using the cyclicity of the trace it is easy to generalize this
result.

\begin{Corollary}
If $\xx=(\x_1,\ldots,\x_{\n-1},0)$,
\begin{equation}\label{Rproduct}
\R^{\tau^j(\xx)}=\G\prod_{i\downarrow_1^{\n-1}}
\psi(\x_i\pw_{1:i+j}\u_{1:i+j-1}\u_{2:i+j}\w_{2:i+j-1}^{-1}),
\end{equation}
where the running indices are considered modulo $\n$.
\end{Corollary}

In what follows we shall use the relations
        \begin{align*}
\G\u_{1:i}&=\u_{1:i}\G,&\G\v_{1:i}&=\v_{1:i}\v_{2:i}
\v_{2:i+1}^{-1}\G\\
\G\u_{2:i}&=\u_{2:i+1}\u_{1:i}\u_{1:i+1}^{-1}\G,&
\G\v_{2:i}&=\v_{2:i+1}\G,
        \end{align*}
 for any $i=1,\ldots,\n\pmod\n$,
which can be deduced from the definition of the operator $\G$.

\subsection{Calculation of the operator $\RR^\xx$}

Introduce more notation:
        \begin{gather}
        \G_{\langle i,j\rangle}\equiv
        \prod_{k\downarrow_i^\n,l\uparrow_1^j}\G_{1:k,2:l} \in
        ((\A^{\otimes\n})^{\otimes\n})^{\otimes2},\notag\\
\label{edefinition}
        \e_{i,j}\equiv\pw_{j:i+j}\u_{j:i+j-1}
        \prod_{k=j}^\n(\u_{k:i+k}^2\u_{k:i+k-1}^{-1}\u_{k:i+k+1}^{-1})
        \in (\A^{\otimes\n})^{\otimes\n},\\
\label{fdefinition}
        \f_{i,j}\equiv\prod_{k=1}^j(\v_{k:i}\v_{k:i+1}^{-1})
        \u_{j:i+1}\w_{j:i}^{-1}\in (\A^{\otimes\n})^{\otimes\n}.
        \end{gather}
        \begin{Proposition}
If $\xx=(\x_1,\ldots,\x_{\n-1},0)$, operator $\RR^\xx$
have the following form
        \begin{equation}\label{rrx}
\RR^\xx=\prod_{i\downarrow_1^\n,j\uparrow_1^\n}\left(
\prod_{k\downarrow_1^{\n-1}}
\psi(\x_k\e_{k,i}\otimes\f_{k,j})\right)\GG,
        \end{equation}
where $\GG\equiv\G_{\langle 1,\n\rangle}$ acts as follows
        \begin{equation}\label{gg}
\GG(f(\ldots,t_{1:i:j},\ldots,t_{2:k:l},\ldots))=
f(\ldots,t_{1:i:j}\prod_{m=1}^\n\frac{t_{2:m:j-i}}{t_{2:m:j-i+1}},
\ldots,t_{2:k:l},\ldots).
        \end{equation}
        \end{Proposition}
\begin{proof} First, calculate
\begin{multline*}
\G_{\langle i,j\rangle}\pw_{1:i:k+i}\u_{1:i:k+i-1}
=\G_{\langle i+1,j\rangle}
\prod_{l\uparrow_1^j}\G_{1:i,2:l}
\pw_{1:i:k+i}\u_{1:i:k+i-1}\\
=\G_{\langle i+1,j\rangle}
\pw_{1:i:k+i}\u_{1:i:k+i-1}
\prod_{m=1}^j(\v_{2:m:k+i}\v_{2:m:k+i+1}^{-1})
\prod_{l\uparrow_1^j}\G_{1:i,2:l}
\\
=\pw_{1:i:k+i}\u_{1:i:k+i-1}\prod_{m=1}^j(\v_{2:m:k}\v_{2:m:k+1}^{-1})
\G_{\langle i,j\rangle},
\end{multline*}
and
\begin{multline*}
\G_{\langle i,j\rangle}
\u_{2:j:k+i}\w_{2:j:k+i-1}^{-1}
=
\G_{\langle i,j-1\rangle}
\prod_{l\downarrow_i^\n}\G_{1:l,2:j}
\u_{2:j:k+i}\w_{2:j:k+i-1}^{-1}
\\
=
\G_{\langle i,j-1\rangle}
\u_{2:j:k+1}\w_{2:j:k}^{-1}
\prod_{m=i}^\n(\u_{1:m:k+m}^2\u_{1:m:k+m-1}^{-1}\u_{1:m:k+m+1}^{-1})
\prod_{l\downarrow_i^\n}\G_{1:l,2:j}
\\
=\u_{2:j:k+1}\w_{2:j:k}^{-1}
\prod_{m=i}^\n(\u_{1:m:k+m}^2\u_{1:m:k+m-1}^{-1}\u_{1:m:k+m+1}^{-1})
\G_{\langle i,j\rangle}.
\end{multline*}
Now we have
\begin{multline*}
\RR^\xx=\prod_{i\downarrow_1^\n,j\uparrow_1^\n}\left(\G_{1:i,2:j}
\prod_{k\downarrow_1^{\n-1}}\psi(\x_k\pw_{1:i:k+i}
\u_{1:i:k+i-1}\u_{2:j:k+i}
\w_{2:j:k+i-1}^{-1})\right)\\
=\prod_{i\downarrow_1^\n,j\uparrow_1^\n}\left(
\prod_{k\downarrow_1^{\n-1}}
\psi(\G_{\langle i,j\rangle}
\x_k\pw_{1:i:k+i}
\u_{1:i:k+i-1}\u_{2:j:k+i}\w_{2:j:k+i-1}^{-1}
\G_{\langle i,j\rangle}^{-1})\right)\GG\\
=\prod_{i\downarrow_1^\n,j\uparrow_1^\n}\left(\prod_{k\downarrow_1^{\n-1}}
\psi(\x_k\e_{k,i}\otimes\f_{k,j})\right)\GG.
\end{multline*}
\end{proof}
\begin{Proposition}
Operators defined by
eqns~(\ref{edefinition}),~(\ref{fdefinition}) satisfy the following
commutation relations
        \begin{gather}\label{permrel}
\e_{i+1,k}\e_{i,k}=\q\e_{i,k}\e_{i+1,k},\quad
\f_{i+1,k}\f_{i,k}=\q\f_{i,k}\f_{i+1,k};\\
\e_{i,k}\e_{j,l}=\q^{(1\pm3)/2}\e_{j,l}\e_{i,k},\quad
\f_{i,l}\f_{j,k}=\q^{(1\pm3)/2}\f_{j,k}\f_{i,l},\quad
k<l,\ |i-j|=(1\mp1)/2;\notag\\
\e_{i,k}\f_{j,k}=\q^{\pm2}\f_{j,k}\e_{i,k},\quad j=k+i-(1\mp1)/2;
        \end{gather}
all other pairs being commutative, and the quadratic constraints
\begin{equation}\label{quadrel}
\e_{i,k+1}\f_{j,k+1}=\f_{j,k}\e_{i,k},\quad k=j-i \pmod\n.
\end{equation}
\end{Proposition}
\begin{proof} It is a straightforward verification.
\end{proof}

\subsection{Realization of a Hopf algebra}

Let $\a_{ij}\equiv2\delta_{ij}-\delta_{|i-j|,1}$ be the
$\mathfrak{sl}_\n$ Cartan matrix. Following \cite{Drinfeld,Jimbo},
consider a Hopf algebra $\Uu$ over $\C(\q)$, generated by elements
$\{\K_i,\L_i,\E_i,\F_i\}_{1\le i<\n}$ subject to the following
relations
        \begin{gather*}
        [\K_i,\K_j]=[\L_i,\L_j]=[\K_i,\L_j]=0,\\
        \K_i\E_j=\q^{\a_{ij}}\E_j\K_i,\quad
        \K_i\F_j=\q^{-\a_{ij}}\F_j\K_i,\\
        \L_i\E_j=\q^{-\a_{ij}}\E_j\L_i,\quad
        \L_i\F_j=\q^{\a_{ij}}\F_j\L_i,        \\
        [\E_i,\F_j]=\delta_{ij}(1-\q^2)(\K_i-\L_i),\\
        \E_i\E_j^2+\E_j^2\E_i=(\q+\q^{-1})\E_j\E_i\E_j\qquad\text{if }
        |i-j|=1,\\
        \F_i\F_j^2+\F_j^2\F_i=(\q+\q^{-1})\F_j\F_i\F_j\qquad\text{if }
        |i-j|=1,\\
        [\E_i,\E_j]=[\F_i,\F_j]=0\qquad\text{if }|i-j|>1.
        \end{gather*}
Note that we use an unusual normalization for the generators.
The coproducts are
        \begin{gather*}
\Delta(\K_i)=\K_i\otimes\K_i,\quad
\Delta(\L_i)=\L_i\otimes\L_i,\\
\Delta(\E_i)=\K_i\otimes\E_i+\E_i\otimes 1,\quad
\Delta(\F_i)=1\otimes\F_i+\F_i\otimes \L_i.
        \end{gather*}
Central elements $\Z_i\equiv\K_i\L_i$ generate a Hopf subalgebra.
Factorization w.r.t. relations $\Z_i=1,\ i=1,\ldots,\n-1$
gives the quantized universal enveloping algebra
$\U_\q(\mathfrak{sl}_\n)$. We construct
a representation of the algebra $\Uu$
in $\C(\q)[t,t^{-1}]^{\otimes\n^2}$, where $\Z_i\ne 1$.

\begin{Theorem}
The following mapping of the generating elements:
        \begin{equation}\label{representation}
        \E_i\mapsto\sum_{j=1}^\n\e_{i,j},\quad
        \F_i\mapsto\sum_{j=1}^\n\f_{i,j}, \quad
        \K_i\mapsto\e_{i,1}\f_{i,1},\quad
        \L_i\mapsto\f_{i,\n}\e_{i,\n}
        \end{equation}
extends to an algebra homomorphism $\kappa\colon\Uu\to\A^{\otimes\n^2}$
\end{Theorem}
\begin{proof}
It is a straightforward checking of the defining relations of the algebra
$\Uu$ by using relations (\ref{permrel})--(\ref{quadrel}).
\end{proof}
\begin{Conjecture}
 The homomorphism $\kappa$ defined by
eqns~(\ref{representation}) is faithful.
\end{Conjecture}
According to the results of \cite{Ber}, the restriction of our 
homomorphism is a particular case of Feigin's homomorphism, and it
is faithful.

\subsection{The universal $R$-matrix}
Certain completion $\hat{\Uu}$ of the algebra $\Uu$
is the Drinfeld double of its' Hopf subalgebra,
so it admits a universal $R$-matrix
\cite{Drinfeld}. To describe the latter
define the $q$-analogues of root vectors
$\{\E_{ij},\F_{ij}\}_{1\le i<j\le\n}$:
        \begin{gather*}
        \E_{i,i+1}\equiv\E_i,\quad
        \F_{i,i+1}\equiv\F_i, \\
        \E_{i,j+1}\equiv\frac{\E_j\E_{i,j}-\q\E_{i,j}\E_j}{1-\q^2},\quad
        \F_{i,j+1}\equiv\frac{\F_j\F_{i,j}-\q\F_{i,j}\F_j}{1-\q^2},\quad i<j,
        \end{gather*}
and Cartan elements $\H_i$ and $\tilde\H_i$ through the equations
        \[
        \K_i=\q^{\H_i},\quad \L_i=\q^{\tilde\H_i}.
        \]
The universal $R$-matrix then has the form \cite{Rosso}:
        \begin{equation}\label{Runiversal}
        \Rr=\prod_{i\downarrow_1^{\n-1}}\prod_{j\downarrow_{i+1}^\n}
        \psi(\E_{ij}\otimes\F_{ij})\q^\T,
        \end{equation}
where
        \[
        \T=\sum_{i,j=1}^{\n-1}\Bar a_{ij}\H_i\otimes\tilde\H_j
        \]
with $\Bar a_{ij}$ being the inverse $\mathfrak{sl}_\n$
Cartan matrix. Element~(\ref{Runiversal}) makes sense
only
for highest weight finite dimensional representations of the
algebra $\Uu$.
The representation defined by eqns~(\ref{representation})
is neither finite dimensional nor highest weight, so
the universal $R$-matrix in the form (\ref{Runiversal}) is not
well defined object. Nevertheless, there is a slight modification
of the definition, when the image of the universal $R$-matrix
makes sense in the $\kappa$-representation.
This modification follows.

Let $\rho_\xx$, where $\xx=(\x_1,\ldots,\x_{\n-1},0)$,
be the automorphism of $\Uu$ defined on the generators
by the formulas:
        \[
        \K_i\mapsto\K_i,\quad \L_i\mapsto\L_i,\quad
        \E_i\mapsto\x_i\E_i,\quad\F_i\mapsto\x_i^{-1}\F_i.
        \]
Note that this automorphism is inner one. Indeed, we have
        \[
        \rho_\xx(.)=\xi_\xx(.)\xi_\xx^{-1},\quad
        \xi_\xx=\prod_{i,j=1}^{\n-1}\x_i^{\Bar a_{ij}\H_j}.
        \]
Denote
        \begin{equation}\label{runivx}
        \Rr^\xx\equiv(\rho_\xx\otimes\id)(\Rr).
        \end{equation}
This element satisfies the multi spectral parameter dependent YBE.
But the spectral parameters here are fictitious,
since they can be removed by a similarity transformation. Their role
is auxiliary: we can look at $\Rr^\xx$ as a generating function
for the coefficients of the power series expansion in
$\x_1,\ldots,\x_{\n-1}$.
\begin{Conjecture}
The following identity holds true
        \begin{equation}\label{mfor}
        \prod_{i\downarrow_1^\n,j\uparrow_1^\n}
        \prod_{k\downarrow_1^{\n-1}}
        \psi(\x_k\e_{k,i}\otimes\f_{k,j})=
        \prod_{i\downarrow_1^{\n-1}}\prod_{j\downarrow_{i+1}^\n}
        \psi(\X_{i,j}\kappa(\E_{i,j})\otimes\kappa(\F_{i,j})),
        \end{equation}
where
        \[
        \X_{i,j}\equiv \prod_{k=i}^{j-1}\x_k.
        \]
\end{Conjecture}
We have verified this identity for few small $\n$.
\begin{Theorem}
Assume that formula~(\ref{mfor}) is true.
Then formula $\kappa\otimes\kappa(\Rr^\xx)=\RR^\xx$
holds as a formal power series equality.
\end{Theorem}
\begin{proof}
Equating formula (\ref{rrx}) to (\ref{runivx}) with substitutions
(\ref{representation}), and using the fact that
$\GG=\q^\T$, we arrive at eqn~(\ref{mfor}).
\end{proof}
\begin{ack}
We would like to thank L.D. Faddeev 
for stimulating discussions and encouragement in this work.
\end{ack}

\end{document}